\documentclass{ws-p8-50x6-00}

\makeatletter
\def\ps@plain{%
     \let\@mkboth\@gobbletwo
     \let\@oddhead\@empty
     \def\@oddfoot{\hfill\thepage\hfill}%
     \let\@evenhead\@empty
     \let\@evenfoot\@oddfoot}
\makeatother

\usepackage{amsmath,amssymb,amscd,amsthm,cite}
\usepackage[mathscr]{eucal}

\newtheorem{theorem}{Theorem}[section]
\newtheorem{proposition}[theorem]{Proposition}
\newtheorem{corollary}[theorem]{Corollary}
\newtheorem{lemma}[theorem]{Lemma}
\newtheorem{definition}{Definition}[section]

\theoremstyle{definition}

\newcommand{\R}{\mathbb{R}}
\newcommand{\Z}{\mathbb{Z}}
\newcommand{\N}{\mathbb{N}}
\newcommand{\Q}{\mathbb{Q}}
\newcommand{\tint}{{\textstyle\int}}
\newcommand{\p}{\partial}

\newcommand{\half}{\tfrac{1}{2}}
\newcommand{\CP}{\mathbb{CP}}
\renewcommand{\*}{\cdot}
\newcommand{\D}{\mathsf{E}}
\newcommand{\E}{\mathcal{E}}
\newcommand{\PP}{\mathsf{P}}
\newcommand{\<}{\langle}
\renewcommand{\>}{\rangle}
\DeclareMathOperator{\Res}{Res}
\newcommand{\eps}{\varepsilon}

\renewcommand{\]}{{]\!]}}

\newcommand{\Om}{\Omega}

\newcommand{\DDD}{\mathcal{D}}
\newcommand{\Mbar}{\overline{\mathcal{M}}}
\renewcommand{\)}{{)\!)}}
\renewcommand{\(}{{(\!(}}
\newcommand{\HH}{\mathsf{K}}
\newcommand{\CL}{\mathcal{L}}

\newcommand{\virt}{{\textup{virt}}}
\DeclareMathOperator{\ev}{ev}

\newcommand{\LL}{\mathcal{L}}
\newcommand{\MM}{\mathcal{M}}

\DeclareMathOperator{\Der}{Der}
\DeclareMathOperator{\ad}{ad}
\newcommand{\bull}{\bullet}

\renewcommand{\AA}{\mathcal{A}}
\newcommand{\AAA}{\mathcal{R}}
\newcommand{\pp}{\delta}
\newcommand{\tpp}{\tilde{\delta}}
\renewcommand{\H}{\mathcal{H}}
\DeclareMathOperator{\DELTA}{[2]}
\newcommand{\ts}{\tilde{s}}
\newcommand{\X}{\mathcal{X}}
\newcommand{\CC}{\mathsf{C}}

\renewcommand{\L}{\mathsf{L}}
\newcommand{\M}{\mathsf{M}}

\begin{document}

\title{The Toda conjecture}

\author{Ezra Getzler}

\address{RIMS, Kyoto University, Kitashirakawa Oiwake-cho, Sakyo-ku, Kyoto
602, Japan}

\address{Northwestern University, Evanston, IL 60208, USA}

\maketitle

Consider the \textbf{Gromov-Witten potential}
$$
F = \sum_{g=0}^\infty \eps^{2g} F_g
$$
of $\CP^1$. Eguchi and Yang\cite{EY} have conjectured that
$$
Z=\exp(\eps^{-2} F )
$$
is a $\tau$-function of the Toda hierarchy. (Similar ideas were also
proposed by Dubrovin, cf.\ \cite{Dubrovin}.) In this paper, we will
explore this conjecture using the bihamiltonian method in the theory
of integrable systems.

Let $P$ be the puncture operator, with descendents $\tau_{k,P}$, $k\ge0$,
and let $Q$ be the operator Poincar\'e dual to a point, with descendents
$\tau_{k,Q}$; denote the corresponding coordinates on the large phase space
by $s_k$ and $t_k$ respectively.

Let $\p$ and $\p_Q$ be the vector fields of differentiation with respect to
$s_0$ and $t_0$, let $\D=e^{\eps\p}$, and introduce the operators
$\nabla=\eps^{-1}(\D^{1/2}-\D^{-1/2})$ and $\DELTA=\D^{1/2}+\D^{-1/2}$. Let
$u$ and $v$ be the functions $\nabla^2F$ and $\nabla\p_QF$.

The \textbf{Toda conjecture} consists of the Toda equation
\begin{equation}
\p_Q^2F = qe^u ,
\end{equation}
whose implications have been studied by Pandharipande\cite{P}, and the
recursion
\begin{equation} \label{r0}
(v\nabla + \DELTA\p_Q) \<\<\tau_{k-1,Q}\>\> = (k+1) \nabla
\<\<\tau_{k,Q}\>\> .
\end{equation}

The large phase space of $\CP^1$ may be identified with the jet-space of
the space with coordinates $u$ and $v$, that is, it has coordinates
$\{\p^nu,\p^nv\}_{n\ge0}$. The Toda conjecture implies that, in these
coordinates, the flows $\p_{k,Q}$ are the flows of the Toda lattice
hierarchy\footnote{The original paper of Eguchi and Yang\cite{EY} had the
incorrect identifications $u=\p\nabla F$ and $v=\p\p_QF$; the corrected
form of the conjecture is found in Eguchi, Hori and Yang\cite{EHY}.}.

Eguchi and Yang also give a matrix integral representation of the
Gromov-Witten potential of $\CP^1$. Studying this representation, Eguchi,
Hori and Yang\cite{EHY} were led to conjecture that $Z$ satisfies a
sequence of constraints $z_n=0$, $n\ge-1$, where $z_{-1}=0$ is the string
equation and $z_0=0$ is Hori's equation. This conjecture is analogous to
the formulation of Witten's KdV conjecture for topological gravity in terms
of an action of the Virasoro algebra, and is now called the
\textbf{Virasoro conjecture} for $\CP^1$; it has recently been proved by
Givental\cite{Givental} (along with its generalization to
higher-dimensional projective spaces).

In Section~\ref{Virasoro}, we discuss the relationship between the Toda
conjecture and the Virasoro conjecture for $\CP^1$; the main result is
that, if \eqref{r0} holds, then the Virasoro conjecture is equivalent of
the following recursion:
\begin{equation} \label{r1}
( v \nabla + \DELTA \p_Q ) \<\<\tau_{k-1,P}\>\> = k \nabla
\<\<\tau_{k,P}\>\> + 2 \nabla \<\<\tau_{k-1,Q}\>\> .
\end{equation}
We show that, when written in terms of the coordinates
$\{\p^nu,\p^nv\}_{n\ge0}$, the commuting flows $\p_{k,P}$ associated to the
descendents of the puncture operator $P$ are Hamiltonian flows; in this
way, we obtain a new hierarchy of Hamiltonians in involution with each
other and with the flows of the Toda lattice.

In Section 6, we show that, in the presence of the Virasoro conjecture, the
Toda conjecture follows once it is known to hold along the submanifold
$\{s_k=0\}_{k>1}$ of the large phase space. Since Okounkov and
Pandharipande have recently proved the Toda conjecture on this submanifold
\cite{OPnew}, the Toda conjecture is established.

Dubrovin and Zhang\cite{DZ} have proved that for homogenous spaces,
and in particular for $\CP^1$, the Virasoro conjecture determines the
Gromov-Witten potential. (They actually prove the analogous result in
the more general context of semisimple Frobenius manifolds.) This
poses the interesting problem of understanding how the Toda conjecture
might follow directly from the Virasoro conjecture.

In this paper, we work over the field $\Q_\eps=\Q\(\eps\)$ of Laurent
polynomials with rational coefficients. The parameter $\eps$ is known
in physics as the loop expansion parameter: this simply means that
integrals over moduli spaces of genus $g$ are weighted by a factor of
$\eps^{2g}$. In the theory of the Toda lattice, $\eps$ is the lattice
spacing --- it is the identification of the lattice spacing with the
genus expansion parameter that lies at the heart of the Toda conjecture.

\numberwithin{equation}{section}

\section{Witten's conjecture}
The Toda conjecture is the analogue for $\CP^1$ of a famous conjecture
of Witten\cite{Witten}, proved by Kontsevich\cite{K}, that the
Gromov-Witten potential of a point is a $\tau$-function of the KdV
hierarchy. (See also Itzykson and Zuber\cite{IZ} and Looijenga\cite{L}
for illuminating discussions of the proof, and Okounkov and
Pandharipande\cite{OP} for an enumerative proof.) In this section, we
recall Witten's conjecture.

Let $\Mbar_{g,n}$ be the moduli space of $n$-pointed stable curves of
arithmetic genus $g$, introduced by Deligne, Mumford and Knudsen; it is
an orbifold of dimension $3(g-1)+n$.

Let $\CL_i$ be the line bundle on $\Mbar_{g,n}$ whose fibre at a
stable curve $(C,z_1,\dots,z_n)$ is the cotangent line $T^*_{z_i}C$,
and let $\psi_i=c_1(\CL_i)\in H^2(\Mbar_{g,n},\Z)$ be its first Chern
class. Witten's conjecture is a formula for the values of the
intersection numbers
$$
\< \tau_{k_1} \dots \tau_{k_n} \>_g = \int_{\Mbar_{g,n}} \psi_1^{k_1} \dots
\psi_n^{k_n} .
$$

It is convenient to assemble these numbers into generating functions
on the large phase space; this is a space with coordinates $t_k$,
$k\ge0$. Denote by $\p_k$ the vector field $\p/\p t_k$ on the large
phase space. The vector field $\p=\p_0$ plays a special role.

Introduce generating functions
$$
F_g = \sum_{n=0}^\infty \frac{1}{n!} \sum_{k_1,\dots,k_n} t_{k_1}
\dots t_{k_n} \< \tau_{k_1} \dots \tau_{k_n} \>_g ,
$$
with partial derivatives $\<\<\tau_{k_1} \dots \tau_{k_n} \>\>_g = \p_{k_1}
\dots \p_{k_n} F_g$, and let $F$ be the total potential
$$
F = \sum_{g=0}^\infty \eps^{2g} F_g ,
$$
with partial derivatives $\<\<\tau_{k_1} \dots \tau_{k_n} \>\> = \p_{k_1}
\dots \p_{k_n} F$.

The total potential $F$ satisfies the \textbf{string equation}
$\CL_{-1}F+\half t_0^2=0$ and the equation $\CL_0F+\frac{1}{8}\eps^2=0$,
where $\CL_{-1}$ and $\CL_0$ are the vector fields
$$
\CL_{-1} = \sum_{k=0}^\infty t_{k+1} \p_k - \p_0 , \quad
\CL_0 = \sum_{k=0}^\infty (k+\half) t_k \p_k - \tfrac{3}{2} \p_1 .
$$
\begin{proposition} \label{basic}
If $f$ is a function on the large phase space such that $\p f$ and
$\CL_{-1}f$ are constant, then $f$ is constant. If in addition,
$\lambda\CL_0f=f$ for some constant $\lambda$, then $f=0$.
\end{proposition}
\begin{proof}
We give an outline of the proof (see Section 3 of Getzler\cite{warsaw} for
more details): if
$$
(\p+\CL_{-1})f = \sum_{k=0}^\infty t_{k+1} \p_kf
$$
is a constant, it follows that $f$ is a constant. Hence $\CL_0f=0$; if
in addition $\lambda\CL_0f=f$, then $f=0$.
\end{proof}

\begin{theorem} \label{jet}
Let $u=\p^2F$. The functions $\p^ku$ $(=\p^{k+2}F)$, $k\ge0$, form a
coordinate system on the large phase space.
\end{theorem}
\begin{proof}
The string equation, in conjunction with the genus $0$ topological
recursion relation
$$
\<\<\tau_k\tau_\ell\tau_m\>\>_0 = \<\<\tau_{k-1}\tau_0\>\>_0
\<\<\tau_0\tau_\ell\tau_m\>\>_0,
$$
implies that $\p_k(\p^\ell u)|_{t_*=0,\eps=0} = \delta_{k\ell}$.
\end{proof}

In the coordinate system $\{\p^ku\}_{k\ge0}$, the vector fields $\CL_{-1}$
and $\CL_0$ have the formulas
\begin{equation} \label{L-1L0}
\CL_{-1} = - \frac{\p}{\p u} , \quad
\CL_0 = - \sum_{k=0}^\infty (\half k+1) \, \p^ku \frac{\p}{\p(\p^ku)} .
\end{equation}

We now recall the definition of the Kortweg-deVries (KdV) hierarchy; this
is a sequence of commuting vector fields on the jet-space of the affine
line. Let $\AA$ be a differential ring (a commutative ring with
differential $\p$), and let $\Psi(\AA)$ be the algebra of
pseudodifferential operators defined over $\AA$: this is the algebra
$$
\Psi(\AA) = \bigcup_{N=0}^\infty \biggl\{ \sum_{i=-\infty}^N a_i \p^i
\biggm| a_i\in\AA \biggr\} ,
$$
with product determined by the relations $\p^i\cdot\p^j=\p^{i+j}$ and
$$
\p^i \cdot a = \sum_{j=0}^\infty \tbinom{i}{j} \p^ja \cdot \p^{i-j} .
$$
Let $A\mapsto A_+$ be the projection on the space of pseudodifferential
operators
$$
\biggl( \sum_{i=-\infty}^N a_i \p^i \biggr)_+ = \sum_{i=0}^N a_i \p^i ,
$$
and let $A_-=A-A_+$.

We are interested in the case where $\AA$ is the algebra of differential
polynomials
$$
\AA = \Q_\eps\[u,\p u,\p^2 u,\dots\] ;
$$
the differential $\p$ acts on the generators as $\p(\p^ku)=\p^{k+1}u$.

The Lax operator is the differential operator
$L=\half\eps^2\p^2+u\in\Psi(\AA)$. There is a unique square root
$D\in\Psi(\AA)$ of $2\eps^{-2}L=\p^2+2\eps^{-2}u$, which commutes with
$L$ and has the form
$$
D = \p + \eps^{-2} u \p^{-1} + \dots .
$$

Let $k$ be a natural number. The \textbf{Lax equation} is the equation
$$
\pp_kL = [(L^kD)_+,L] ,
$$
or equivalently, $\pp_kL=- [(L^kD)_-,L]$. From these two equations,
we see that $\pp_kL$ is an element of $\AA$. Write
$$
L^kD = \sum_{i=-\infty}^{2k+1} a_i(k) \p^i ,
$$
where $a_i(k)$ is an element of $\AA$. The differential polynomial
$f_k=\eps^2 a_{-1}(k)$is called the $k$th \textbf{Gelfand-Dickii
polynomial}. Since $\pp_kL$ is the constant term in the commutator
$-[(L^kD)_-,L]$, we see that
$$
\pp_kL=\p f_k\in\AA .
$$

The Lax equation determines a derivation of $\AA$, defined on generators by
$$
\pp_k(\p^nu) = \p^n(\pp_kL) = \p^{n+1}f_k .
$$
The essential property of the Lax equation is that these flows commute:
since $\pp_kD=[(L^kD)_+,D]$, it follows that
$\pp_m(L^nD)_+=[(L^mD)_+,L^nD]_+$, and we see that
\begin{multline*}
[\pp_m,\pp_n]L = \pp_m[(L^nD)_+,L] - \pp_n[(L^mD)_+,L] \\
\begin{aligned}
{}&= [\pp_m(L^nD)_+,L] + [(L^nD)_+,\pp_mL] - [\pp_n(L^mD)_+,L] -
[(L^mD)_+,\pp_nL] \\
{}&= [[(L^mD)_+,L^nD]_+,L] - [[(L^nD)_+,L^mD]_+,L] - [[(L^mD)_+,(L^nD)_+],L] \\
{}&= [[L^mD,L^nD]_+,L] = 0 .
\end{aligned}
\end{multline*}

\begin{theorem} \label{GD}
Let $\HH = \tfrac{1}{8} \eps^2 \p^3 + u\p + \half \p u$. The
differential polynomials $f_k$ are characterized by two properties:
the recursion $\HH f_{k-1} = \p f_k$ holds, and $f_k$ has vanishing
constant term.
\end{theorem}
\begin{proof}
It is clear that the recursion $\HH f_{k-1} = \p f_k$ determines $f_k$
up to a constant, since the kernel of the linear map $\p:\AA\to\AA$
consists of multiples of the identity. Furthermore, $f_k$ has
vanishing constant term, since $D|_{u=0}=\p$, hence $(L^kD)_-|_{u=0}$
vanishes. It remains to prove the recursion.

The vanishing of the coefficient of $\p^i$ in the equation $[L,L^kD]=0$
gives
$$
\p a_{i-1}(k) = - \half \p^2a_i(k) + \eps^{-2} \sum_{j=1}^\infty
\tbinom{i+j}{j} a_{i+j}(k)\,\p^ju .
$$
Taking $i=-2$ and $i=-1$, we see that
$$
\p a_{-3}(k) = - \half \p^2 a_{-2}(k) - \eps^{-2} \p u \, a_{-1}(k) =
\bigl( \tfrac14 \p^3 - \eps^{-2} \p u \bigr) a_{-1}(k) .
$$

By considering the coefficient of $\p^i$ in the equations
$L^{k+1}D=(L^kD)L$, we see that
$$
a_i(k+1) = \half \eps^2 a_{i-2}(k) + \sum_{j=0}^\infty \tbinom{i+j}{j}
\p^ju\,a_{i+j}(k) ,
$$
and in particular, that $a_{-1}(k+1) = \half \eps^2 a_{-3}(k) + u \,
a_{-1}(k)$. Taking a derivative of this equation gives
$$
\p a_{-1}(k+1) = \half \eps^2 \p a_{-3}(k) + \p u \, a_{-1}(k) + u \, \p
a_{-1}(k) ,
$$
and the recursion follows.
\end{proof}

Since $f_0=u$, we see from Theorem \ref{GD} that
$f_1=\tfrac{1}{8}\eps^2\p^2u + \tfrac{3}{4}u^2$. In particular, $\pp_0=\p$,
while $\pp_1u=\tfrac{1}{8}(\eps^2\p^3u+12u\p u)$ is the KdV equation.

Let
$$
\alpha_k = \<\<\tau_0\tau_k\>\> - \tfrac{2^k}{(2k+1)!!} f_k ;
$$
in particular, $\alpha_0=0$. Witten's conjecture has a number of equivalent
formulations:
\begin{gather}
\framebox[21pc]
{for all $k\ge0$, $\HH\<\<\tau_0\tau_k\>\> = (k+\half) \,
\p\<\<\tau_0\tau_{k+1}\>\>$} \tag{i} \\
\framebox[21pc]
{for all $k\ge0$, $\alpha_k=0$} \tag{ii} \\
\framebox[21pc]
{for all $k\ge0$, the vector field $\p_k$ equals
$\frac{2^k}{(2k+1)!!}\pp_k$}
\tag{iii}
\end{gather}
It is obvious that (ii) implies (i); let us show that (i) implies
(ii). Since $\CL_{-1}u=-1$, we see that
$\CL_{-1}(L^kD)=-(k+\half)L^{k-1}D$, hence $\CL_{-1}f_k = - (k+\half)
\, f_{k-1}$; it follows that
$\CL_{-1}\alpha_k=-\alpha_{k-1}$. Likewise, we may prove by induction,
using the recursion $\HH f_{k-1}=\p f_k$, that
$$
\CL_0f_k = - \sum_{n=0}^\infty (\half n+1) \, \p^nu
\tfrac{\p}{\p(\p^nu)} f_k = -(k+1) f_k ,
$$
hence that $\CL_0\alpha_k=-(k+1)\alpha_k$. Together, (i) and Theorem
\ref{GD} imply that $\HH\alpha_{k-1}=(k+\half)\p\alpha_k$. Using
Proposition \ref{basic}, we may now argue by induction that $\alpha_k=0$.

Likewise, it is obvious that (ii) implies (iii), since
$\p\alpha_k=\bigl(\p_k-\tfrac{2^k}{(2k+1)!!}\delta_k\bigr)u$. The
proof of the converse is similar to the proof that (i) implies (ii).

The involutivity of the vector fields $\pp_k$ is essential to the
formulation of Witten's conjecture: it is seen to be an integrability
condition for the existence of the coordinates $t_k$ on the large phase
space. Clearly, this conjecture determines $\p F$; in combination with the
equations $\CL_{-1}F+\half t_0^2=0$ and $\CL_0F+\tfrac18\eps^2=0$, it
follows from Proposition \ref{basic} that it determines $F$.

\section{The Toda lattice}

In this section, we introduce the Toda lattice, in the form in which it
enters into the Toda conjecture: the limit in which the lattice spacing
$\eps$ is infinitesimal (Takasaki and Takebe\cite{TT}). We follow the
approach of Kupershmidt\cite{K}.

If $(\AA,\p)$ is a commutative algebra with derivation $\p$ over
$\Q_{\eps,q}=\Q_{\eps,q}$, let $\D:\AA\to\AA$ be the automorphism
$e^{\eps\p}$. Let $\Phi(\AA)$ be the algebra of twisted Laurent series
with coefficients $\AA$; as a vector space, $\Phi(\AA)$ is the space
of Laurent series $\AA\(\Lambda^{-1}\)$, and the product is given by
the formula
$$
\sum_i a_i \Lambda^i \* \sum_j b_j \Lambda^j = \sum_{i,j}
(\D^{-j/2}a_i)(\D^{i/2}b_j) \Lambda^{i+j} .
$$
Extend the automorphism $\D$ to $\Phi(\AA)$ by letting $\D$ act
trivially on $\Lambda$. For example, if $\AA=C^\infty_c(\R)$ and
$\p=d/dt$, then $\Phi(\AA)$ is a continuum limit of the algebra of
infinite matrices $(M_{ij})_{i,j\in\Z}$ such that $M_{ij}=0$ for
$|i-j|\gg0$.

Let $A\mapsto A_+$ be the projection on $\Phi(\AA)$ defined the
formula
$$
\biggl( \sum_{i=-\infty}^N a_i\Lambda^i \biggr)_+ = \sum_{i=0}^N a_i\Lambda^i ,
$$
and let $A_-=A-A_+$.

The commutative algebra with derivation which we will use is
$$
\AA = \Q_{\eps,q}\[e^u,\p^nu,\p^nv\mid n\ge0\] .
$$
The derivation $\p$ acts on the generators in the evident way:
\begin{align*}
\p e^u &= e^u\,\p u , & \p(\p^nu) &= \p^{n+1}u , & \p(\p^nv)=\p^{n+1}v .
\end{align*}
The kernel of the operators $\p$ and $\nabla$ on $\AA$ equals
$\Q_{\eps,q}$. Let $\PP$ be the infinite-order differential operator
$$
\PP = \frac{\p}{\nabla} = \sum_{g=0}^\infty
\frac{\eps^{2g}(2^{1-2g}-1)B_{2g}}{(2g)!}  \, \p^{2g} = 1 - \tfrac{1}{24}
\, \eps^2 \, \p^2 + O(\eps^4) .
$$
Obviously, $\p=\PP\circ\nabla$.

\begin{definition}
The \textbf{Lax operator} of the Toda lattice is
$$
L = \Lambda + v + qe^u \Lambda^{-1} \in \Phi(\AA) .
$$
\end{definition}

Define elements $p_k(n)\in\AA$, $n\ge0$, $k\in\Z$, as follows:
$$
L^n = \sum_{k=-\infty}^n p_k(n) \Lambda^k .
$$
\begin{lemma}
$p_{-1}(n)=qe^up_1(n)$
\end{lemma}
\begin{proof}
Taking the coefficient of $\Lambda^0$ in the equation $[L,L^n]=0$, we
see that
$$
\nabla p_{-1}(n) = \nabla(qe^up_1(n)) .
$$
Thus $p_{-1}(n)-qe^up_1(n)\in\Q_{\eps,q}$. But when $u=v=0$, the Lax
operator equals $\Lambda+q\Lambda^{-1}$, so that
$p_{-k}(n)|_{u=v=0}=q^kp_k(n)|_{u=v=0}$.
\end{proof}

The following proposition shows that the functions $p_k(n)$ have certain
homogeneity properties.
\begin{proposition} \label{LL}
Let $e$ and $\E$ be the vector fields
\begin{equation} \label{Ee}
e = \frac{\p}{\p v} \text{ and } \E = \sum_{n=0}^\infty \p^nv \,
\frac{\p}{\p(\p^nv)} + 2 \, \frac{\p}{\p u} .
\end{equation}
Then $e(p_k(n))=np_k(n-1)$ and $\E(p_k(n))=(n-k)p_k(n$).
\end{proposition}
\begin{proof}
The vector fields $e$ and $\E$ both commute with $\p$, hence with the
operator $\D$. It follows that they induce derivations of the algebra
$\Phi(\AA)$. Since $e(L)=1$, we see that $e(L^n)=nL^{n-1}$, and expanding
in powers of $\Lambda$, that $e(p_k(n))=np_k(n-1)$.

Likewise, since $\E(L)=v+2qe^u\Lambda^{-1}=(1-\Lambda\p_\Lambda)L$ and
$\Lambda\p_\Lambda$ is a derivation of $\Phi(\AA)$, we see that
$$
\E(L^n)=(n-\Lambda\p_\Lambda)L^n;
$$
expanding in powers of $\Lambda$, it follows that $\E(p_k(n))=(n-k)p_k(n)$.
\end{proof}

The $n$th \textbf{Toda flow} is determined by the Lax equation
$$
\pp_nL = \eps^{-1}[L^n_+,L] = - \eps^{-1}[L^n_-,L] .
$$
Since $L^n_+$ involves only positive powers of $\Lambda$ and $L^n_-$ involves
only negative powers of $\Lambda$, the coefficient of $\Lambda^i$ in
$\pp_nL$ vanishes unless $i$ equals $0$ or $-1$, and
$$
\pp_nL = \nabla p_{-1}(n) + qe^u \nabla p_0(n)\Lambda^{-1} .
$$
There is a unique derivation $\pp_n$ of the algebra $\AA$, which
commutes with $\p$ and is characterized by the formulas $\pp_nu=\nabla
p_0(n)$ and $\pp_nv=\nabla p_{-1}(n)$. In particular, the derivation
$\delta_1$ is the original Toda flow:
$$
\pp_1u = \nabla v , \quad \pp_1v = q \nabla e^u .
$$
Eliminating $v$, we obtain the Toda equation $\pp^2_1u=q\nabla^2e^u$.

Being defined by Lax equations, the Toda flows commute: by the formula
$\pp_mL^n_+=[L^m_+,L^n]_+$, we see that
\begin{align*}
\eps^2 [\pp_m,\pp_n]L &= \pp_m[L^n_+,L] - \pp_n[L^m_+,L] \\
&= [\pp_m(L^n)_+,L] + [L^n_+,\pp_mL] - [\pp_n(L^m)_+,L] -
[L^m_+,\pp_nL] \\
&= [[L^m_+,L^n]_+,L] - [[L^n_+,L^m]_+,L] - [[L^m_+,L^n_+],L] \\
&= [[L^m,L^n]_+,L] = 0 .
\end{align*}

Let $\Om_\AA$ be the algebra of differential forms
$$
\Om^\bull_\AA = \AA[d(\p^nu),d(\p^nv)\mid n\ge0] ,
$$
generated over $\AA$ by Grassmann variables $\{d(\p^nu),d(\p^nv)\}_{n\ge0}$
of degree $1$. There is a unique derivation $\p$ on $\Om_\AA$ which agrees
with the derivation $\p$ on $\AA$ and commutes with the exterior
differential $d$.

The space of \textbf{functional differential forms} is the cokernel of
$\p$:
$$
\AAA^\bull = \Om^\bull_\AA/\p\Om^\bull_\AA .
$$
The image of an element $\alpha\in\Om^\bull_\AA$ in $\AAA^\bull$ is denoted
$\tint\alpha\,dt$; this notation is intended to indicate that integration
by parts is permitted under the integral sign:
$$
\tint\p\alpha\wedge\beta \, dt = - \tint\alpha\wedge\p\beta \, dt .
$$
The exterior differential $d$ on $\Om^\bull_\AA$ induces a
differential on $\AAA^\bull$. Elements of $\AAA^0$ are called
\textbf{functionals}.

There is a natural identification between $\AAA^1$ and
$\AA\,du\oplus\AA\,dv$, since
$$
\text{$\tint f_k\,d(\p^ku) \, dt = \tint (-\p)^kf_k \, du \, dt$ and $\tint
g_k\,d(\p^kv) \, dt = \tint (-\p)^kg_k \, dv \, dt$.}
$$
Under this identification, the exterior differential $d:\AAA^0\to\AAA^1$
may be written
$$
d=\delta_u\,du+\delta_v\,dv ,
$$
where $\delta_u$ and $\delta_v:\AAA^0\to\AA$ are the \textbf{variational
derivatives} (also known as \textbf{Euler-Lagrange operators})
\begin{equation} \label{EL}
\delta_u = \sum_{k=0}^\infty (-\p)^k \frac{\p}{\p(\p^ku)} , \quad
\delta_v = \sum_{k=0}^\infty (-\p)^k \frac{\p}{\p(\p^kv)} .
\end{equation}
\begin{lemma} \label{Res}
The \textup{\textbf{residue}} $\Res:\Phi(\Om_\AA)\to\AAA$, defined by the
formula
$$
\Res(a_n\Lambda^n) = \begin{cases} \tint a_0 \, dt , & n=0 , \\ 0 ,
& \text{otherwise,}
\end{cases}$$
vanishes on graded commutators and satisfies $d\Res(f)=\Res(df)$.
\end{lemma}
\begin{proof}
It is clear that $d\Res(f)=\Res(df)$, and that $\Res [ a\Lambda^k ,
b\Lambda^\ell ]$ vanishes unless $k+\ell=0$, while
$\Res[a\Lambda^k,b\Lambda^{-k}]=(\D^{k/2}-\D^{-k/2})(ab)
\in\p\Om^\bull_\AA$.
\end{proof}

The functionals $h_n = \tfrac{1}{n+1} \Res(L^{n+1}) = \tfrac{1}{n+1}
\tint p_0(n+1) \, dt\in \AAA^0$ are the \textbf{Hamiltonians} of the
Toda lattice hierarchy.
\begin{proposition} \label{H}
$\delta_uh_n=p_{-1}(n)$ and $\delta_vh_n=p_0(n)$
\end{proposition}
\begin{proof}
It follows from Lemma \ref{Res} that $dh_n=\Res(L^ndL)$. Since
$$dL=dv+qe^udu\Lambda^{-1},$$ we see that
$\Res(L^ndL)=\tint\bigl(p_0(n)dv+qe^up_1(n)du\bigr)dt$.
\end{proof}

In the dispersionless limit $\eps\to0$, the algebra $\Phi(\AA)$ becomes the
commutative algebra of Laurent series $\AA\(\Lambda^{-1}\)$, and it is
straightforward to calculate a generating function for these Hamiltonians
(Fairlie and Strachan\cite{FS}).
\begin{proposition}
$$
\lim_{\eps\to\infty} \, \sum_{n=0}^\infty t^n p_0(n) = \bigl( 1 - 2tv +
t^2(v^2-4qe^u) \bigr)^{-1/2}
$$
\end{proposition}
\begin{proof}
Let $w=v-t^{-1}$. Let $\gamma$ be a small circular contour around the
origin of the complex plane. Using the residue formula, we may write the
generating function which we wish to calculate as
{\footnotesize
\begin{multline*}
\lim_{\eps\to\infty} \sum_{n=0}^\infty t^n p_0(n) =
\frac{1}{2\pi i} \int_\gamma \frac{1}{1-t(\Lambda+v+qe^u\Lambda^{-1})}
\frac{d\Lambda}{\Lambda} \\
\begin{aligned}
{}&= - \frac{1}{2\pi i t} \int_\gamma \frac{d\Lambda}
{(\Lambda+\frac12 w+\frac12(w^2-4qe^u)^{1/2})
(\Lambda+\frac12 w-\frac12(w^2-4qe^u)^{1/2})} \\
{}&= \frac{1}{2\pi i t(w^2-4qe^u)^{1/2}} \int_\gamma \Bigl(
\frac{d\Lambda}{\Lambda+\frac12 w+\frac12(w^2-4qe^u)^{1/2}} -
\frac{d\Lambda}{\Lambda+\frac12 w - \frac12(w^2-4qe^u)^{1/2}} \Bigr) .
\end{aligned}
\end{multline*}}
The two poles of the integrand are at $\half w+\half(w^2-4qe^u)^{1/2}=O(t)$
and $\half w-\half(w^2-4qe^u)^{1/2}=-t^{-1}+O(1)$ respectively; thus, for
sufficiently small values of $t$, only the first contributes, with residue
$1$, and the generating function equals
$t^{-1}(w^2-4qe^u)^{-1/2}=(1-2tv+t^2(v^2-4qe^u))^{-1/2}$.
\end{proof}

\begin{corollary}
Let $P_n(x)$ be the $n$th Legendre polynomial. Then
$$
p_0(n) = (v^2-4qe^u)^{n/2} P_n\bigl( v/(v^2-4qe^u)^{1/2} \bigr) .
$$
\end{corollary}
\begin{proof}
The Legendre polynomials have generating function
$$
\sum_{n=0}^\infty t^n P_n(x) = (1-2xt+t^2)^{-1/2} .
$$
Setting $x=v/(v^2-4qe^u)^{1/2}$ and $t=t(v^2-4qe^u)^{1/2}$, the result
follows.
\end{proof}

The following theorem of Kupershmidt\cite{K} gives a pair of formulas for
the derivations $\pp_n$, in terms of $h_n$, and $h_{n-1}$ respectively. (He
also proves a third formula for $\pp_n$, in terms of $h_{n-2}$.)
\begin{theorem} \label{Kupershmidt}
Let $\CC_u=\eps^{-1}q(\D^{1/2}\*e^u\*\D^{1/2} -
\D^{-1/2}\*e^u\*\D^{-1/2})$. Then
$$
\pp_n \begin{bmatrix} v \\ u \end{bmatrix} =
\begin{bmatrix} 0 & \nabla \\ \nabla & 0 \end{bmatrix}
\begin{bmatrix} \delta_vh_n \\ \delta_uh_n \end{bmatrix} =
\begin{bmatrix} \CC_u & v\nabla \\ \nabla v & \eps^{-1}(\D-\D^{-1})
\end{bmatrix}
\begin{bmatrix} \delta_vh_{n-1} \\ \delta_uh_{n-1} \end{bmatrix} .
$$
\end{theorem}
\begin{proof}
We have already proved the first of these formulas. To prove the second
formula, observe that by the equations $L^n=L^{n-1}L=LL^{n-1}$, we have
{\small
\begin{align*}
p_k(n) &= \D^{-1/2}p_{k-1}(n-1) + (\D^{k/2}v)\,p_k(n-1) +
q(\D^{(k+1)/2}e^u)\,\D^{1/2}p_{k+1}(n-1) \\
&= \D^{1/2}p_{k-1}(n-1) + (\D^{-k/2}v)\,p_k(n-1) +
q(\D^{-(k+1)/2}e^u)\,\D^{-1/2}p_{k+1}(n-1) .
\end{align*}}
In particular, the equality of these two formulas for $p_0(n)$ shows
that
\begin{equation} \label{p-1}
\nabla p_{-1}(n) = q\nabla (e^up_1(n)) .
\end{equation}
Further, taking $\D^{1/2}$ times the first formula minus $\D^{-1/2}$ times
the second, with $k$ respectively $-1$ and $0$, gives
\begin{align*}
\nabla p_{-1}(n) &= v\nabla p_{-1}(n-1) + \CC_up_0(n-1) \\
\nabla p_0(n) &= \nabla(vp_0(n-1)) +
\eps^{-1}q(\D-\D^{-1})(e^up_1(n-1)) .
\end{align*}
Since $q(\D-\D^{-1})(e^up_1(n-1))=(\D-\D^{-1})p_{-1}(n-1)$ by Lemma
\ref{p-1}, the result follows.
\end{proof}

\begin{corollary} \label{DDD}
$p_0(n) = vp_0(n-1) + \DELTA p_{-1}(n-1)$
\end{corollary}
\begin{proof}
Since $\nabla p_0(n) = \nabla \bigl( vp_0(n-1) + \DELTA p_{-1}(n-1)
\bigr)$, we see that
$$
p_0(n) - vp_0(n-1) - \DELTA p_{-1}(n-1) \in \Q_{\eps,q} .
$$
To show that this vanishes, we evaluate it at the point $u=v=0$: since
$L|_{u=v=0}=\Lambda+q\Lambda^{-1}$, we see that
$$
p_0(n)|_{u=v=0} = (qp_1(n-1)+p_{-1}(n-1))|_{u=v=0}=2\,p_{-1}(n-1)|_{u=v=0}
,
$$
as required.
\end{proof}

\section{Hamiltonian operators and the Toda lattice}

In this section, we introduce the variational Schouten Lie algebra; this is
an infinite dimensional analogue of the usual Schouten Lie algebra,
developed by Dorfman and Gelfand\cite{GDorfman}. We explain how it may be
used to give an alternative approach to the Toda lattice. For more details,
together with applications to other hierarchies, see
Dorfman\cite{Dorfman}, Getzler\cite{darboux}, Manin\cite{Manin} and
Olver\cite{Olver}.

We start by introducing the free graded commutative algebra
$\Lambda_\infty$ over the algebra $\AA$, with generators
$\{\p^n\theta_v,\p^n\theta_u\}_{n\ge0}$ of degree $1$. The derivation $\p$
is extended to $\Lambda_\infty$ by the formulas
\begin{align*}
\p(\p^n\theta_v) &= \p^{n+1}\theta_v , & \p(\p^n\theta_u) &=
\p^{n+1}\theta_u .
\end{align*}
The kernel of $\p$ is spanned by $1\in\Lambda_\infty$, and the cokernel of
$\p$ is denoted $\LL$. Denote the image of an element $F\in\Lambda_\infty$
in $\LL$ by $\tint F\,dt$.

The variational derivatives $\delta_u$ and $\delta_v$ on $\Lambda_\infty$
are defined by the same formulas \eqref{EL} as on $\AA$, while the
associated Grassmann variational derivatives, which we denote by $\delta^u$
and $\delta^v$, are defined by the formulas
\begin{align*}
\delta^u &= \sum_{n=0}^\infty (-\p)^n \frac{\p}{\p(\p^n\theta_u)} &
\delta^u &= \sum_{n=0}^\infty (-\p)^n \frac{\p}{\p(\p^n\theta_v)} .
\end{align*}
Since all of these operators vanish on the image of $\p$, they descend to
linear operators from $\LL$ to $\Lambda_\infty$. The \textbf{Schouten
bracket} is the bilinear operation on $\LL$ defined by the formula
$$
\bigl[ \tint F\,dt , \tint G\,dt \bigr] = \tint \bigl( \delta^uF \,
\delta_uG + \delta^vF \, \delta_vG + (-1)^{|F|} ( \delta_uF \, \delta^uG +
\delta_vF \, \delta^vG ) \bigr) dt .
$$
With this bracket, the graded vector space $\LL$ becomes a graded Lie
algebra: $[f,g]=(-1)^{|f|\,|g|}[g,f]$ and $[f,[g,h]] = [[f,g],h] +
(-1)^{(|f|+1)(|g|+1)} [g,[f,h]]$. (Proofs of these formulas may be
found in \cite{darboux}.)

The Lie subalgebra $\LL^1$ of $\LL^\bull$ is isomorphic to the Lie algebra
of derivations of $\AA$ commuting with $\p$, under the map
$$
\X=\tint(f\,\theta_u+g\,\theta_v)\,dt \in \LL^1 \mapsto \sum_{n=0}^\infty
\biggl( \p^nf \, \frac{\p~}{\p(\p^nu)} + \p^ng \, \frac{\p~}{\p(\p^nv)}
\biggr) \in \Der(\AA) .
$$
For example, the vector fields $e$ and $\E$ of \eqref{Ee} correspond to
$\tint\theta_v\,dt$ and $\tint(v\theta_v+2\theta_u)\,dt$; we will denote
these elements of $\LL^1$ by $e$ and $\E$ as well.

An element $\H$ of $\LL^2$ defines a graded derivation $\delta_\H$ of
degree $1$ on the graded Lie algebra $\LL$, by the formula $\delta_\H =
[\H,-]$. Let $f\in\LL^0=\AAA^0$ be a functional; the derivation of $\AA$
commuting with $\p$ which corresponds to the element $\delta_\H f$ of
$\LL^1$ is called the \textbf{Hamiltonian vector field} associated to
$f$. The \textbf{Poisson bracket} on the space of functionals
$\LL^0\cong\AAA^0$ is defined by the formula
$$
\{f,g\}_\H = [\delta_\H f,g] .
$$
If $[\H,\H]=0$, $\H$ is called a \textbf{Hamiltonian operator}.
\begin{proposition} \label{Hamiltonian}
If $\H$ is a Hamiltonian operator, $\delta_\H$ is a differential, the
bracket $\{f,g\}_\H$ is a Lie bracket, and $[\delta_\H f,\delta_\H g] =
\delta_\H\{f,g\}_\H$.
\end{proposition}
\begin{proof}
By the graded Jacobi rule,
$$
\delta_\H(\delta_\H f) = [\H,[\H,f]] = \half [[\H,\H],f] = 0 ,
$$
showing that $\delta_\H$ is a differential.

We have
$$
\{f,g\}_\H + \{g,f\}_\H = [\delta_\H f,g] + [\delta_\H g,f] =
\delta_\H[f,g] .
$$
Since $[f,g]$ vanishes, we see that $\{f,g\}_\H$ is antisymmetric.

By the graded Jacobi rule, we have
\begin{multline*}
\{\{f,g\}_\H,h\}_\H = [\delta_\H[\delta_\H f,g],h] = [[\delta_\H f
,\delta_\H g],h] \\ = [\delta_\H f,[\delta_\H g,h]] - [\delta_\H
g,[\delta_\H f,h]] = \{f,\{g,h\}_\H\}_\H - \{g,\{f,h\}_\H\}_\H .
\end{multline*}
Finally, we have $[\delta_\H f,\delta_\H g] = \delta_\H[\delta_\H f,g] =
\delta_\H\{f,g\}_\H$.
\end{proof}

A \textbf{bihamiltonian structure} is a pair of Hamiltonian operators
$(\H,\H_0)$ such that $\H+\lambda\H_0$ is a Hamiltonian operator for all
$\lambda$, or equivalently, such that $[\H,\H_0]=0$. We now reformulate
Theorem \ref{Kupershmidt} in terms of a bihamiltonian structure.
\begin{theorem}
The operators
\begin{align*}
\H &= \half \tint \begin{bmatrix} \theta_v & \theta_u \end{bmatrix}
\begin{bmatrix} \CC_u & v\nabla \\ \nabla v & \eps^{-1}(\D-\D^{-1})
\end{bmatrix} \begin{bmatrix} \theta_v \\ \theta_u \end{bmatrix} \, dt \\
&= \tint \bigl( \eps^{-1} \theta_u \, \D\theta_u + v\theta_v\nabla\theta_u
+ \eps^{-1} qe^u(\D^{-1/2}\theta_v)(\D^{1/2}\theta_v) \bigr) \, dt
\text{, and} \\
\H_0 &= \half \tint \begin{bmatrix} \theta_v & \theta_u \end{bmatrix}
\begin{bmatrix} 0 & \nabla \\ \nabla & 0 \end{bmatrix}
\begin{bmatrix} \theta_v \\ \theta_u \end{bmatrix} \, dt
= \tint \theta_v\nabla\theta_u \, dt
\end{align*}
give a bihamiltonian structure.
\end{theorem}
\begin{proof}
It is clear that $\delta_v\H=\theta_v(\nabla\theta_u)$. Using the formula
$$
\delta^v = \sum_{i=-\infty}^\infty \D^{-i/2}
\frac{\p}{\p(\D^{i/2}\theta_v)} ,
$$
we see that $\delta^v\H=v\nabla\theta_u + \CC_u\theta_v$. It follows
that
\begin{align*}
\delta^v\H \, \delta_v\H &= \CC_u\theta_v \, \theta_v (\nabla\theta_u) \\
&= \eps^{-1}q(\D^{1/2}e^u)(\D\theta_v)\,\theta_v (\nabla\theta_u)
- \eps^{-1}q(\D^{-1/2}e^u)(\D^{-1}\theta_v)\,\theta_v (\nabla\theta_u) .
\end{align*}
Likewise, since $\delta_u\H=\eps^{-1}
qe^u(\D^{-1/2}\theta_v)(\D^{1/2}\theta_v)$ and
$\delta^u\H=\nabla(v\theta_v)+\DELTA\nabla\theta_u$, it follows that
$$
\delta^u\H \, \delta_u\H = \eps^{-1}qe^u \DELTA(\nabla\theta_u)
(\D^{-1/2}\theta_v)(\D^{1/2}\theta_v) .
$$
We conclude that
\begin{align*}
\delta^v\H \, \delta_v\H + \delta^u\H \, \delta_u\H &= \eps^{-1}q \bigl(
(\D-\D^{-1/2})(e^u(\D^{1/2}\theta_v)(\D^{-1/2}\theta_v)\theta_u) \\ &\quad
+ (1-\D^{-1/2})(e^u(\D^{1/2}\theta_v)(\D^{-1/2}\theta_v)(\D\theta_u)) \\
&\quad+
(\D^{1/2}-1)(e^u(\D^{1/2}\theta_v)(\D^{-1/2}\theta_v)(\D^{-1}\theta_u))
\bigr) ,
\end{align*}
and hence that $[\H,\H]=0$. The proof that
$[\H+\lambda\H_0,\H+\lambda\H_0]=0$ is a formal consequence of this
formula, since $\H+\lambda\H_0$ is obtained from $\H$ by translating $v$ by
$\lambda$.
\end{proof}

The Hamiltonian operators $\H_0$ and $\H$ have nice commutation relations
with the vector fields $e$ and $\E$: it is easily checked that
$[e,\H]=\H_0$ and $[e,\H_0]=0$, and that $[\E,\H]=0$ and
$[\E,\H_0]=-\H_0$.

Denote the Poisson bracket associated to $\H$ by $\{f,g\}_\H$, and the
Poisson bracket associated to $\H_0$ by $\{f,g\}_0$; then Theorem
\ref{Kupershmidt} may be reformulated as the identity
\begin{equation} \label{recurse}
[\H_0,h_n] = [\H,h_{n-1}] ,
\end{equation}
or equivalently, as the pair of equations
\begin{align}
\nabla(\delta_vh_n) &= \nabla(v\delta_vh_{n-1}) +
\eps^{-1}(\D-\D^{-1})\delta_uh_{n-1} \label{qq0} \\
\nabla(\delta_uh_n) &= \CC_u\delta_vh_{n-1} + v\nabla\delta_uh_{n-1}
. \label{qq1}
\end{align}
Using \eqref{recurse}, we obtain another proof that the flows associated to
the Hamiltonians $h_n$ commute; this proof has the advantage that it does
not depend on the Lax equation, and so applies to prove the involutivity of
more general hierarchies.
\begin{proposition} \label{Magri}
Let $f_n$ and $g_n$ be sequences of Hamiltonians such that
$[\H_0,f_n]=[\H,f_{n-1}]$ and $[\H_0,g_n]=[\H,g_{n-1}]$ for $n>0$, and
$\{f_0,f\}=0$ for all $f$. Then $\{f_m,g_n\}_0 = 0$ for all
$m,n\ge0$. In particular, the Hamiltonian vector fields associated to
$f_m$ and $g_n$ commute.
\end{proposition}
\begin{proof}
If $m>0$, we have
\begin{align*}
\{f_m,g_n\}_0 &= [[\H_0,f_m],g_n] = [[\H,f_{m-1}],g_n] = -
[[\H,g_n],f_{m-1}] \\ &= - [[\H_0,g_{n+1}],f_{m-1}]
= \{f_{m-1},g_{n+1}\}_0 .
\end{align*}
Thus, it suffices to prove the proposition for $m=0$, for which it is
clear.
\end{proof}

It is an immediate consequence of Proposition \ref{LL} that
\begin{equation} \label{Ehk}
\text{$[e,h_k]=kh_{k-1}$ and $[\E,h_k]=(k+1)h_k$.}
\end{equation}

If $a\in\Z$, let $\LL^i(a)$ be the generalized eigenspace
$$
\LL^i(a) = \bigcup_{n>0} \ker\bigl( (\ad(\E)+i-a-1)^n \bigr) ;
$$
then $[\LL(a),\LL(b)] = \LL(a+b)$. Since $[\E,\H_0]=-\H_0$, we see that
$\H_0\in\LL(0)$, hence the differential $\delta_0$ preserves the graded
subspace $\LL(a)$. Also, we see that $\H\in\LL(1)$.
\begin{theorem} \label{main}
The complex $\LL(a)$ has vanishing cohomology unless $a=-1$ or $0$, while
the nonzero cohomology groups of $\LL(-1)$ and $\LL(0)$ are as follows:
\begin{align*}
H^0(\LL(-1),\delta_0) &= \<\tint 1\,dt , \tint u\,dt \> &
H^0(\LL(0),\delta_0) &= \< \tint v\,dt \> \\
H^1(\LL(-1),\delta_0) &= \< \tint \theta_v\,dt \> &
H^1(\LL(0),\delta_0) &= \< \tint \theta_u\,dt , \tint
(u\theta_u-v\theta_v)\,dt \> \\
& & H^2(\LL(0),\delta_0) &= \< \tint \theta_u\theta_v\,dt \>
\end{align*}
\end{theorem}
\begin{proof}
Write elements of the cone $\MM$ of the linear map
$\nabla:\Lambda_\infty/\Q_{\eps,q}\to\Lambda_\infty$ as $(F,G)$, where
$F\in\Lambda_\infty$ and $G\in\Lambda_\infty/\Q_{\eps,q}$. If $D$ is an
linear operator from $\Lambda_\infty/\Q_{\eps,q}$ to $\Lambda_\infty$, let
$\iota D:\MM\to\MM$ be the operator
$$
\iota D(F,G) = (DG,0) .
$$
The differential of $\MM$ equals $\iota\nabla$.

The differential $\delta_0$ on $\LL$ lifts to a differential
$$
\delta = \sum_{n=0}^\infty \bigl( (\nabla\p^n\theta_u)\p_{\p^nv} +
(\nabla\p^n\theta_v)\p_{\p^nu} \bigr)
$$
on $\MM$, and the map $\tau:\MM\to\LL$ which sends $(F,G)$ to $\tint F\,dt$
is a chain homotopy equivalence between the complexes
$(\MM,\delta+\iota\nabla)$ and $(\LL,\delta_0)$.

Let $S$ be the chain homotopy
$$
S = \sum_{n=0}^\infty \bigl( (\PP\p^nv)\p_{\p^{n+1}\theta_u} +
(\PP\p^nu)\p_{\p^{n+1}\theta_v} \bigr) ,
$$
let $U=\iota(v\p_{\theta_u}+u\p_{\theta_v})$, and let $P$ be the semisimple
operator
$$
P = \sum_{n=0}^\infty \bigl( (\p^nu)\p_{\p^nu} + (\p^nv)\p_{\p^nv} \bigr) +
\sum_{n=1}^\infty \bigl( (\p^n\theta_u)\p_{\p^n\theta_u} +
(\p^n\theta_v)\p_{\p^n\theta_v} \bigr) .
$$
We have $[\delta+\iota\nabla,S]=e^{-U}Pe^U=P+U$. Thus, the cohomology of
the complex $(\MM,d+\iota\nabla)$ equals the kernel
$$
\< (1,0), (\theta_u,0), (\theta_v,0), (\theta_u\theta_v,0), (u,-\theta_v),
(v,-\theta_u), (u\theta_u-v\theta_v,\theta_u\theta_v) \> \subset \MM
$$
of $[\delta+\iota\nabla,S]$. Applying $\tau$, the theorem follows.
\end{proof}

\section{The Toda conjecture}

Having introduced the Toda lattice hierarchy in the last two sections, we
can now formulate the Toda conjecture. Recall the definition of the
Gromov-Witten invariants of $\CP^1$. (A good review of the subject is Manin
\cite{frobenius}.) Let $\Mbar_{g,n,d}$ be the moduli stack of stable maps
of genus $g$ and degree $d$, with $n$ marked points, to $\CP^1$. Let
$\ev_i:\Mbar_{g,n,d}\to\CP^1$, $1\le i\le n$, be the map
$$
\ev_i(f:C\to\CP^1,z_1,\dots,z_n) = f(z_i) .
$$
defined by evaluating a stable map $f:C\to\CP^1$ at the $i$th marked
point $z_i\in C$. Let
$$
[\Mbar_{g,n,d}]^\virt \in H_{2(2g-2+2d+n)}(\Mbar_{g,n,d},\Q)
$$
be the virtual fundamental class.

\begin{definition}
Let $\CL_i$ be the line bundle on $\Mbar_{g,n,d}$ whose fibre at the
stable map $(f:C\to\CP^1,z_1,\dots,z_n)$ is the line $T^*_{z_i}C$, and
let $\psi_i=c_1(\CL_i)$ be its first Chern class.
\end{definition}

Let $\gamma_P\in H^0(\CP^1,\Z)$ and $\gamma_Q\in H^2(\CP^1,\Z)$ be the
cohomology classes Poincar\'e dual to the fundamental class and to a
point. Given $k_i\in\N$ and $a_i\in\{P,Q\}$, define
$$
\<\tau_{k_1,a_1}\dots\tau_{k_n,a_n}\>_g = \sum_{d=0}^\infty q^d
\int_{[\Mbar_{g,n,d}]^\virt} \ev_1^*\gamma_{a_1}\dots\ev_n^*\gamma_{a_n}
\cup \psi_1^{k_1} \dots \psi_n^{k_n} \in \Q[q] .
$$
We write $P$ and $Q$ instead of $\tau_{0,P}$ and $\tau_{0,Q}$.

The \textbf{large phase space} $\M$ is the formal manifold with coordinates
$\{s_k,t_k\}_{k\ge0}$. Define
$$
t_k^a = \begin{cases} s_k , & a=P , \\ t_k , & a=Q . \end{cases}
$$
The genus $g$ \textbf{Gromov-Witten potential} $F_g$ of $\CP^1$ is the
generating function on the large phase space given by the formula
$$
F_g = \sum_{n=0}^\infty \frac{1}{n!}  \sum_{\substack{k_1,\dots,k_n \\
a_1,\dots,a_n}} \prod_{i=1}^n t_{k_i}^{a_i} \int_{[\Mbar_{g,n,d}]^\virt}
\<\tau_{k_1,a_1}\dots\tau_{k_n,a_n}\>_g ,
$$
and $F= \sum_{g=0}^\infty \eps^{2g} F_g$ is the total Gromov-Witten
potential.

Denote the constant vector fields $\p/\p t_k^a$ on the large phase space by
$\p_{k,a}$; in particular, write $\p$ and $\p_Q$ for $\p_{0,P}$ and
$\p_{0,Q}$. Just as in the theory of the Toda lattice, denote the
operator $e^{\eps\p}$ by $\D$, and the operators
$\eps^{-1}(\D^{1/2}-\D^{-1/2})$ and $\D^{1/2}+\D^{-1/2}$ by $\nabla$ and
$\DELTA$. The partial derivatives of $F$ are denoted
$\<\<\tau_{k_1,a_1}\dots\tau_{k_n,a_n}\>\> =
\p_{k_1,a_1}\dots\p_{k_n,a_n}F$.

The potential $F$ satisfies the \textbf{string equation}
$\CL_{-1}F+s_0t_0=0$ (Witten\cite{Witten}) and \textbf{Hori's equation}
$\CL_0F+s_0^2=0$ (Hori\cite{Hori}), where $\CL_{-1}$ and $\CL_0$ are the
vector fields
\begin{align}
\CL_{-1} &= \sum_{k=0}^\infty \bigl( t_{k+1} \p_{k,Q} + s_{k+1} \p_{k,P}
\bigr) - \p \label{L-1} \\
\CL_0 &= \sum_{k=0}^\infty \bigl( (k+1) \, t_k \p_{k,Q} + k \, s_k \p_{k,P}
+ 2 \, s_{k+1} \p_{k,Q} \bigr) - \p_{1,P} - 2\,\p_Q . \label{L0}
\end{align}
The analogue of Proposition \ref{basic} holds, with essentially the same
proof.
\begin{proposition} \label{Basic}
If $f$ is a function on the large phase space such that $\nabla f$ and
$\CL_{-1}f$ are constant, then $f$ is constant. If in addition,
$\lambda\CL_0f=f$ for some constant $\lambda$, then $f=0$.
\end{proposition}

There is also an analogue for $\CP^1$ of Theorem \ref{jet}.
\begin{theorem} \label{Jet}
Let $u=\nabla^2F$ and $v=\nabla\p_QF$. The functions
$\{\p^nu,\p^nv\}_{n\ge0}$ form a coordinate system on the large phase
space. The origin $s_k=t_k=0$ of the large phase space has coordinates
\begin{align*}
\p^nv &= \begin{cases} 1 , & n=1 , \\ 0 , & \text{otherwise,}
\end{cases} &
\p^nu &= 0 .
\end{align*}
\end{theorem}
\begin{proof}
The string equation, in conjunction with the genus $0$ topological
recursion relation
$$
\<\<\tau_{k,a}\tau_{\ell,b}\tau_{m,c}\>\>_0 = \<\<\tau_{k-1,a}P\>\>_0
\<\<Q\tau_{\ell,b}\tau_{m,c}\>\>_0 + \<\<\tau_{k-1,a}Q\>\>_0
\<\<P\tau_{\ell,b}\tau_{m,c}\>\>_0 ,
$$
implies that $\p_{k,P}(\p^nu)|_{s_*=t_*=0} = \p_{k,Q}(\p^nv)|_{s_*=t_*=0} =
\delta_{kn} + O(\eps)$.

At the origin of the large phase space, we have
$$
\p^nv = \sum_{k=1}^\infty \sum_{g=0}^\infty
\frac{\eps^{2g+2k-1}}{2^{2k-2}(2k-1)!} \<P^{2k-1+n}Q\>_g .
$$
The moduli space $\Mbar_{g,2k+n,d}$ only contributes to this sum if it has
virtual dimension $1$, that is, if $2(g-1+d+k)+n=1$. The only solution of
this equation is $\Mbar_{0,3,0}\cong\CP^1$, which contributes the
coefficient $1$ to $\p v$. The argument for $\p^nu$ is similar: at the
origin in the large phase space,
$$
\p^nu = \sum_{k=1}^\infty \sum_{g=0}^\infty
\frac{\eps^{2g+2k-2}}{2^{2k-1}(2k)!} \<P^{2k+n}\>_g .
$$
The only contributions to this sum come from moduli spaces
$\Mbar_{g,2k+n,d}$ such that $2(g-1+d+k)+n=0$; this equation has no
solutions.
\end{proof}

The vector fields $\CL_{-1}$ and $\CL_0$ commute with $\p$; written in the
coordinate system $\{\p^nu,\p^nv\}_{n\ge0}$, they may be identified with
$-e$ and $-\E$ (see \eqref{Ee}).

Introduce the constraints
\begin{align*}
\alpha_{k,Q}^v &= \nabla\<\<\tau_{k-1,Q}\>\> - \tfrac{1}{k!} \delta_vh_k , &
\alpha_{k,Q}^u &= \p_Q\<\<\tau_{k-1,Q}\>\> - \tfrac{1}{k!} \delta_uh_k .
\end{align*}
Observe that $\alpha_{1,Q}^v$ and $\alpha_{0,Q}^u$ vanish. By Proposition
\ref{LL}, we have
\begin{align} \label{La}
\CL_{-1}\alpha^*_{k,Q} &= -\alpha^*_{k-1,Q} , & \CL_0\alpha_{k,Q}^v &=
-k\alpha_{k,Q}^v , & \CL_0\alpha_{k,Q}^u &= -(k+1)\alpha_{k,Q}^u .
\end{align}

The Toda conjecture describes the vector fields $\p_{k,Q}$ in the
coordinate system $\{\p^nu,\p^nv\}_{n\ge0}$. It may be formulated as any
one of the equivalent conditions in the following theorem.
\begin{theorem} \label{Y}
Let $\DDD$ be
the differential operator
$$
\DDD = v\nabla + \DELTA\,\p_Q .
$$
The following are equivalent:
\begin{align}
&\framebox[24pc]{for all $k>0$, $\DDD\<\<\tau_{k-1,Q}\>\> = (k+1) \,
\nabla \<\<\tau_{k,Q}\>\>$, and $\p_Q^2F=qe^u$} \tag{i} \\
&\framebox[24pc]{for all $k>0$,
$k!\,\nabla\<\<\tau_{k-1,Q}\>\>=\delta_vh_k$ and
$k!\,\p_Q\<\<\tau_{k-1,Q}\>\>=\delta_uh_k$} \tag{ii} \\
&\framebox[24pc]{for all $k>0$, the vector field $k!\,\p_{k-1,Q}$ equals
$\pp_k$} \tag{iii}
\end{align}
\end{theorem}
\begin{proof}
Introduce the constraint
$$
y_k = \DDD \<\<\tau_{k,Q}\>\> - (k+2) \nabla \<\<\tau_{k+1,Q}\>\> .
$$
We may reformulate the conditions of the theorem in terms of the
constraints $y_k$, $\alpha^u_{k,Q}$ and $\alpha^v_{k,Q}$ as follows:
\begin{gather}
\framebox[15pc]{for all $k\ge0$, $y_k=0$, and $\alpha_{1,Q}^u=0$} \tag{i} \\
\framebox[15pc]{for all $k>0$, $\alpha_{k,Q}^*=0$} \tag{ii} \\
\framebox[15pc]{for all $k>0$, $\nabla\alpha_{k,Q}^*=0$} \tag{iii}
\end{gather}
The third of these reformulations follows from the formulas
$k!\,\nabla\alpha_{k,Q}^v=(k!\,\p_{k-1,Q}-\pp_k)u$ and
$k!\,\nabla\alpha_{k,Q}^u=(k!\,\p_{k-1,Q}-\pp_k)v$.

\begin{lemma} \label{y}
We have $y_{k-1}=v\alpha_{k,Q}^v + \DELTA\alpha_{k,Q}^u - (k+1)
\alpha_{k+1,Q}^v$ and
$$
\p_Qy_{k-1} = \CC_u\alpha_{k,Q}^v + v\nabla\alpha_{k,Q}^u -
(k+1)\nabla\alpha_{k+1,Q}^u .
$$
\end{lemma}
\begin{proof}
The first formula follows from the calculation
\begin{align*}
(k+1) \alpha_{k+1,Q}^v &= (k+1) \nabla\<\<\tau_{k,Q}\>\> - \tfrac{1}{k!}
\delta_vh_{k+1} \\ &= \DDD\<\<\tau_{k-1,Q}\>\> - \tfrac{1}{k!} (
v\delta_vh_k + \DELTA\delta_uh_k ) - y_{k-1} .
\end{align*}
To prove the second formula, observe that
\begin{align*}
(k+1)\nabla\alpha_{k+1,Q}^u &= (k+1)\nabla\p_Q\<\<\tau_{k,Q}\>\> -
\tfrac{1}{k!} \nabla\delta_uh_{k+1} \\
&= \p_Q\DDD\<\<\tau_{k-1,Q}\>\> - \tfrac{1}{k!} \CC_u\delta_vh_k -
\tfrac{1}{k!} v\nabla\delta_uh_k - \p_Qy_{k-1} .
\end{align*}
Since
$\p^2_Q\<\<\tau_{k-1,Q}\>\>=\p_{k-1,Q}e^u=e^u\nabla^2\<\<\tau_{k-1,Q}\>\>$
and $\p_Qv=\nabla e^u$, we see that
\begin{align*}
\p_Q\DDD\<\<\tau_{k-1,Q}\>\> &= (\p_Qv)\nabla\<\<\tau_{k-1,Q}\>\> +
v\p_Q\nabla\<\<\tau_{k-1,Q}\>\> + \DELTA\p^2_Q\<\<\tau_{k-1,Q}\>\> \\
&= (\nabla e^u) \nabla\<\<\tau_{k-1,Q}\>\> +
v\p_Q\nabla\<\<\tau_{k-1,Q}\>\> + \DELTA ( e^u \nabla^2\<\<\tau_{k-1,Q}\>\>
) .
\end{align*}
A short calculation shows that $(\nabla e^u)f + \DELTA(e^u\nabla f)
=\CC_uf$, proving the second formula.
\end{proof}

By \eqref{La} and Proposition \ref{Basic}, we see that if
$\nabla\alpha_{i,Q}^*=0$ for $i\le k$, then $\alpha_{k,Q}^*=0$; in
particular, (iii) implies (ii). By Lemma \ref{y}, it is clear that (ii)
implies (i) and that (i) implies (iii).
\end{proof}

Condition (iii) of Theorem \ref{Y} has recently been proved by
Okounkov and Pandharipande \cite{OPnew} on the submanifold
$\{s_k=0\}_{k>1}$ of the large phase space; as we will see in Section
6, in conjunction with the Virasoro conjecture, this establishes the
Toda conjecture.

\section{The Toda conjecture and the Virasoro conjecture} \label{Virasoro}

The Virasoro conjecture for $\CP^1$ says that the functions $z_k$,
$k\ge-1$, vanish, where $z_k$ is given by the formula
\begin{multline*}
z_k = \sum_{m=0}^\infty \bigl( c_k^{m+1} t_m \<\<\tau_{m+k,Q}\>\> + c_k^m
\ts_m \<\<\tau_{m+k,P}\>\> + 2 d_k^m \ts_m \<\<\tau_{m+k-1,Q}\>\> \bigr) \\
+ \sum_{i+j=k} i!\,j! \bigl( \eps^2 \<\<\tau_{i-1,Q}\tau_{j-1,Q}\>\> +
\<\<\tau_{i-1,Q}\>\> \<\<\tau_{j-1,Q}\>\> \bigr) + \delta_{k,-1}s_0t_0 +
\delta_{k,0} s_0^2 .
\end{multline*}
Here, $c_k^m=m(m+1)\dots(m+k)$, and $d_k^m=e_k(m,\dots,m+k)$, where
$e_k(x_0,\dots,x_k)$ is the $k$th elementary symmetric function, and
$\ts_m=s_m-\delta_{m,1}$. This conjecture was made by Eguchi, Hori and
Yang\cite{EHY}, motivated by their matrix integral representation of the
Gromov-Witten potential of $\CP^1$. The string equation and Hori's equation
are the special cases with $k=-1$ and $k=0$ respectively. The formulas
\begin{align} \label{Lz}
\CL_{-1}z_k &= -(k+1)z_{k-1} , & \CL_0z_k &= -kz_k ,
\end{align}
may be proved by direct calculation.

The Virasoro conjecture for $\CP^1$ has been proved by Givental
\cite{Givental}; his proof uses Kontsevich's localization theorem for
Gromov-Witten invariants of toric varieties \cite{Kon}, together with
results from Dubrovin's theory of semisimple Frobenius manifolds
\cite{Dubrovin}.

In  this section, we study the constraints
$$
x_k = \DDD \<\<\tau_{k,P}\>\> - (k+1) \nabla \<\<\tau_{k+1,P}\>\> - 2
\nabla \<\<\tau_{k,Q}\>\> .
$$
Our main result is that if the Toda and Virasoro conjectures hold,
then the constraints $x_k$ vanish.
\begin{theorem} \label{master}
\begin{align*}
\DDD z_k - \nabla z_{k+1} &= \sum_{m=0}^\infty \bigl( c_k^{m+1} \, t_m \,
y_{m+k} + c_k^m \, \ts_m \, x_{m+k} + 2 \, d_k^m \, \ts_m \, y_{m+k-1}
\bigr) \\ &+ \sum_{i+j=k} i! \, j!  \bigl( \eps^2 \p_{i-1,Q} +
\DELTA\<\<\tau_{i-1,Q}\>\> \bigr) y_{j-1}
\end{align*}
\end{theorem}
\begin{proof}
We make use of the following formulas:
$[\p_{i,Q},\DDD]=\nabla\<\<\tau_{i,Q}Q\>\>\nabla$,
\begin{gather*}
[\DDD,t_k] = \delta_{k,0}\DELTA , \quad
[\DDD,s_k] = \half\delta_{k,0}\bigl(\eps^2\nabla\p_Q+v\DELTA\bigr) , \\
\nabla(fg) = \half\DELTA f\,\nabla g + \half \nabla f \, \DELTA g , \\
\DDD(fg) = \half \DELTA f\,\DDD g + \half \DDD f \, \DELTA g + \half
\eps^2 \p_Q (\nabla f\,\nabla g) .
\end{gather*}
It follows that
{\small
\begin{align*}
& \DDD\sum_{m=0}^\infty ( c_k^{m+1}t_m\<\<\tau_{m+k,Q}\>\> +
c_k^m\ts_m\<\<\tau_{m+k,P}\>\> + 2d_k^m\ts_m\<\<\tau_{m+k-1,Q}\>\> ) \\
{}&= \sum_{m=0}^\infty ( c_k^{m+1}t_m\DDD\<\<\tau_{m+k,Q}\>\>
+ c_k^m\ts_m\DDD\<\<\tau_{m+k,P}\>\> + 2d_k^m\ts_m\DDD\<\<\tau_{m+k-1,Q}\>\>
) \\ 
{}&\quad+ c_k^1 \DELTA \<\<\tau_{k,Q}\>\> +
d_k^0 (\eps^2\nabla\p_Q+v\DELTA)
\<\<\tau_{k-1,Q}\>\> \\
{}&= \sum_{m=0}^\infty (
 c_{k+1}^{m+1}t_m\nabla\<\<\tau_{m+k+1,Q}\>\>
+ c_{k+1}^m\ts_m\nabla\<\<\tau_{m+k+1,P}\>\> +
2d_{k+1}^m\ts_m\nabla\<\<\tau_{m+k,Q}\>\> ) \\
{}&\quad+ \sum_{m=0}^\infty ( c_k^{m+1}t_my_{m+k}
+ c_k^m\ts_mx_{m+k} + 2d_k^m\ts_my_{m+k-1} ) \\
{}&\quad+ c_k^1 \DELTA\<\<\tau_{k,Q}\>\> + d_k^0
(\eps^2\nabla\p_Q+v\DELTA) \<\<\tau_{k-1,Q}\>\> \\
{}&= \nabla \sum_{m=0}^\infty ( c_{k+1}^{m+1}t_m\<\<\tau_{m+k+1,Q}\>\> +
c_{k+1}^m\ts_m\<\<\tau_{m+k+1,P}\>\> + 2d_{k+1}^m\ts_m\<\<\tau_{m+k,Q}\>\>
) \\
{}&\quad+ \sum_{m=0}^\infty ( c_k^{m+1}t_my_{m+k}
+ c_k^m\ts_mx_{m+k} + 2d_k^m\ts_my_{m+k-1} ) \\
&\quad+ d_k^0 (\eps^2\nabla\p_Q+v\DELTA)\<\<\tau_{k-1,Q}\>\>
\end{align*}
\begin{align*}
\DDD \sum_{i+j=k} & i! \, j! ( \eps^2 \<\<\tau_{i-1,Q}\tau_{j-1,Q}\>\> +
\<\<\tau_{i-1,Q}\>\> \<\<\tau_{j-1,Q}\>\> ) \\
{}&= \sum_{i+j=k} i! \, j! \bigl( \eps^2 \p_{i-1,Q} \DDD
\<\<\tau_{j-1,Q}\>\> - \eps^2 [\p_{i-1,Q},\DDD] \<\<\tau_{j-1,Q}\>\> \\
{}& \qquad + \DELTA\<\<\tau_{i-1,Q}\>\>
\DDD\<\<\tau_{j-1,Q}\>\> + \half \eps^2 \p_Q (\nabla\<\<\tau_{i-1,Q}\>\>
\nabla\<\<\tau_{j-1,Q}\>\>) \bigr) \\
&= \sum_{i+j=k} i! \, (j+1)! ( \eps^2
\nabla\<\<\tau_{i-1,Q}\tau_{j,Q}\>\> + \DELTA\<\<\tau_{i-1,Q}\>\>
\nabla\<\<\tau_{j,Q}\>\> ) \\
&\quad+ \sum_{i+j=k} i! \, j! ( \eps^2 \p_{i-1,Q} +
\DELTA\<\<\tau_{i-1,Q}\>\> ) y_{j-1} \\
{}&= \nabla \sum_{i+j=k+1} i! \, j! ( \eps^2
\<\<\tau_{i-1,Q}\tau_{j-1,Q}\>\> + \<\<\tau_{i-1,Q}\>\>\<\<\tau_{j-1,Q}\>\>
) \\
&\quad- k! ( \eps^2 \nabla\p_Q + v\DELTA )
\<\<\tau_{k-1,Q}\>\> + \sum_{i+j=k} i! \, j! ( \eps^2 \p_{i-1,Q} +
\DELTA\<\<\tau_{i-1,Q}\>\> ) y_{j-1} .
\end{align*}}
The theorem follows for $k>0$ on adding the results of these two
calculations; the cases $k=-1$ and $k=0$ are similar, and we leave them to
the reader.
\end{proof}

\begin{corollary}
If $y_k=0$ for all $k\ge0$ and $z_k=0$ for all $k\ge-1$, then $x_k=0$ for
all $k\ge0$.
\end{corollary}
\begin{proof}
If $y_k$ and $z_k$ vanish for all $k$, then the formula of Theorem
\ref{master} becomes, for $k\ge-1$,
$$
x_{k+1} = \sum_{m=1}^\infty \tbinom{k+m}{m-1} \, s_m \, x_{k+m} .
$$
The result follows by induction on the order of vanishing of the
constraints $x_k$ at the origin of the large phase space.
\end{proof}

Assuming the Toda and Virasoro conjectures, we will now show that there are
Hamiltonians $g_k\in\AAA^0$ such that
$\p_{k,P}u=\tfrac{1}{k!}\nabla\delta_vg_k$ and
$\p_{k,P}v=\tfrac{1}{k!}\nabla\delta_ug_k$. We may construct the
Hamiltonian $g_0$ explicitly: the equations $\p u=\nabla\delta_vg_0$ and
$\p v=\nabla\delta_ug_0$ have the solution
\begin{equation} \label{g0}
g_0 = \tint u \PP v \, dt = \sum_{g=0}^\infty
\frac{\eps^{2g}(2^{1-2g}-1)B_{2g}}{(2g)!}  \tint u \, \p^{2g}v \, dt .
\end{equation}
Note that $\delta_0g_0$ lies in the centre of $\LL$, since
\begin{equation} \label{centre}
[\delta_0g_0,\tint F\,dt] = [\tint(\p u\,\theta_u+\p
v\,\theta_v)\,dt,\tint F\,dt] = \tint \p F\,dt = 0 .
\end{equation}
We have not been able to find an explicit formula for the Hamiltonians
$g_k$; instead, we construct them using a method due to Gelfand and Dorfman
\cite{GDorfman}.
\begin{theorem}
There is a unique sequence of Hamiltonians $g_k\in\LL^0(k)$ starting with
$g_0=\tint u\PP v\,dt$ such that
$[\H_0,g_k]=[\H,g_{k-1}-\frac{2}{k}h_{k-1}]$.
\end{theorem}
\begin{proof}
Uniqueness is clear: if $k>0$, the recursion
$[\H_0,g_k]=[\H,g_{k-1}-\frac{2}{k}h_{k-1}]$ determines $g_k$ up to an
element of $H^0(\LL(k),\delta_0)$, and this cohomology group vanishes by
Theorem \ref{main}.

We will construct $g_k$ by induction.
\begin{lemma}
The vector field $\X_k=[\H,g_{k-1}-\frac{2}{k}h_{k-1}]\in\LL^1(k)$
satisfies $\delta_0\X_k=0$.
\end{lemma}
\begin{proof}
We have
$$
\delta_0\X_k = [\H_0,[\H,g_{k-1}-\tfrac{2}{k}h_{k-1}]] = -
[\H,[\H_0,g_{k-1}-\tfrac{2}{k}h_{k-1}]] .
$$
If $k>1$, we have
$$
[\H,[\H_0,g_{k-1}-\tfrac{2}{k}h_{k-1}]] =
[\H,[\H,g_{k-2}-2(\tfrac{1}{k}+\tfrac{1}{k-1})h_{k-2}]] = 0 .
$$
If $k=1$, we see that
$\delta_0\X_1=-[\H,[\H_0,g_0-2h_0]]=[\delta_0g_0,\H]$, which vanishes by
\eqref{centre}.
\end{proof}

It follows from this lemma that the cohomology class of $\X_k$ is an
element of $H^1(\LL(k),\delta_0)$. For $k>0$, this cohomology group
vanishes by Theorem~\ref{main}, hence there is an element $g_k$ of
$\LL^0(k)$ such that $\delta_0g_k=\X_k$.
\end{proof}

\begin{corollary} \label{magri}
Denote the Hamiltonian vector field $\delta_0g_k$ associated to the
Hamiltonian $g_k$ by $\tpp_k$:
\begin{align*}
\tpp_kv &= \nabla\delta_ug_k , & \tpp_ku &= \nabla\delta_vg_k .
\end{align*}
The flows $\tpp_k$ commute with each other, and with the flows $\pp_k$.
\end{corollary}
\begin{proof}
Apply Theorem \ref{Magri} to the sequences of Hamiltonians $h_k$ and $g_k$.
\end{proof}

It is not hard to find an explicit formula for $g_1$:
\begin{equation} \label{g1}
g_1 = \tint\bigl( \half u\PP(v^2+\DELTA e^u) + \half v(\DELTA\PP-2)v -
2qe^u \bigr) \, dt .
\end{equation}
The equation
$$
[\H_0,g_1] = [\H,g_0-2h_0] = [\H,g_0] - 2[\H_0,h_1]
$$
amounts by \eqref{qq0} and \eqref{qq1} to the pair of equations
\begin{align*}
\nabla(\delta_vg_1) &= \nabla(v\delta_vg_0+\DELTA\delta_ug_0-2\delta_vh_1)
= \nabla(v\PP u+\DELTA(\PP-2)v) \\
\nabla(\delta_ug_1) &= \CC_u\delta_vg_1 +
v\nabla\delta_ug_1-2\nabla\delta_uh_1 \\
&= \half q\nabla(e^u\DELTA\PP u) + \half q\DELTA(e^u\p u) + v\nabla\PP v -
2q\nabla e^u \\
&= q\nabla \bigl( \half e^u\DELTA\PP u + \half \PP(v^2+\DELTA e^u) - 2 e^u
\bigr) ,
\end{align*}
and it is easily seen that these are satisfied.

\begin{proposition}
For $k>0$, $[e,g_k]=kg_{k-1}$. For $k\ge0$, $$[\E,g_k]=(k+1)g_k+2h_k.$$
\end{proposition}
\begin{proof}
We prove these formulas by induction: it is easily checked, using the
explicit formula \eqref{g1} for $g_1$, that $[e,g_1]=g_0$, and, using the
explicit formula \eqref{g0} for $g_0$, that $[\E,g_0]=g_0+2h_0$.

We have
\begin{align*}
\delta_0[e,g_k] &= [e,\delta_0g_k] = [e,[\H,g_{k-1}-\tfrac{2}{k}h_{k-1}]] \\
&= [[e,\H],g_{k-1}-\tfrac{2}{k}h_{k-1}] +
[\H,[e,g_{k-1}-\tfrac{2}{k}h_{k-1}]] \\
&= \delta_0(g_{k-1}-\tfrac{2}{k}h_{k-1}) +
[\H,(k-1)g_{k-2}-\tfrac{2(k-1)}{k}h_{k-2}] \\
&= \delta_0(g_{k-1}-\tfrac{2}{k}h_{k-1}) +
\delta_0((k-1)g_{k-1}+\tfrac{2}{k}h_{k-1}) = \delta_0(kg_{k-1}) .
\end{align*}
Since the cohomology group $H^0(\LL(k-1),\delta_0)$ vanishes for $k>1$, the
formula for $[e,g_k]$ follows.

The argument for $\E$ is similar: we have
\begin{align*}
\delta_0[\E,g_k] &= [\delta_0\E,g_k] + [\E,\delta_0g_k] =
\delta_0g_k + [\E,[\H,g_{k-1}-\tfrac{2}{k}h_{k-1}]] \\
&= \delta_0g_k + [\H,[\E,g_{k-1}-\tfrac{2}{k}h_{k-1}]] \\
&= \delta_0g_k + [\H,kg_{k-1}] = \delta_0((k+1)g_k+2h_k) .
\end{align*}
Since the cohomology group $H^0(\LL(k),\delta_0)$ vanishes for $k>0$, the
formula for $[\E,g_k]$ follows.
\end{proof}

Introduce the constraints
\begin{align*}
\alpha_{k,P}^v &= \nabla\<\<\tau_{k,P}\>\> - \tfrac{1}{k!} \delta_vg_k , &
\alpha_{k,P}^u &= \p_Q\<\<\tau_{k,P}\>\> - \tfrac{1}{k!} \delta_ug_k .
\end{align*}
By the definitions of $u$ and $v$, we see that $\alpha_{0,P}^v$ and
$\alpha_{0,P}^u$ vanish.
\begin{lemma} \label{Lb}
For $k>0$, $\CL_{-1}\alpha^*_{k,P}=-\alpha^*_{k-1,P}$ and
\begin{align*}
\CL_0\alpha_{k,P}^v &= -k\alpha_{k,P}^v - 2\alpha_{k,Q}^v , &
\CL_0\alpha_{k,P}^u &= -(k+1)\alpha_{k,P}^u - 2\alpha_{k,Q}^u .
\end{align*}
\end{lemma}
\begin{proof}
By the string equation,
$\CL_{-1}\<\<\tau_{k,P}\>\>=-\<\<\tau_{k-1,P}\>\>$. Since
$[e,\delta_v]=[e,\delta_u]=0$, we see that
$\CL_{-1}\delta_vg_k=-k\delta_vg_{k-1}$ and
$\CL_{-1}\delta_ug_k=-k\delta_ug_{k-1}$. This shows that
$\CL_{-1}\alpha^*_{k,P}=-\alpha^*_{k,P}$.

By Hori's equation $z_0=0$, we see that
\begin{align*}
\CL_0\nabla\<\<\tau_{k,P}\>\> &= -k\nabla\<\<\tau_{k,P}\>\> -
2\nabla\<\<\tau_{k-1,Q}\>\> , \\
\CL_0\p_Q\<\<\tau_{k,P}\>\> &= -(k+1)\p_Q\<\<\tau_{k,P}\>\>
-2\p_Q\<\<\tau_{k-1,Q}\>\> .
\end{align*}
Since $[\E,\delta_v]=-\delta_v$ and $[\E,\delta_u]=0$, we see that
$\CL_{-1}\delta_vg_k=-k\delta_vg_k-2\delta_vh_k$ and
$\CL_0\delta_ug_k=-(k+1)\delta_ug_k-2\delta_uh_k$. This yields the formulas
for $\CL_0\alpha^*_{k,P}$.
\end{proof}

\begin{theorem} \label{X}
Assume that the Toda conjecture holds. Each of the following conditions are
equivalent to the Virasoro conjecture:
\begin{align}
& \framebox[23pc]{for all $k>0$, $\DDD\<\<\tau_{k-1,P}\>\> =
k\nabla\<\<\tau_{k,P}\>\> + 2\<\<\tau_{k-1,Q}\>\>$} \tag{i} \\
& \framebox[23pc]{for all $k>0$,
$k!\,\nabla\<\<\tau_{k,P}\>\>=\delta_vg_k$ and
$k!\,\p_Q\<\<\tau_{k,P}\>\>=\delta_ug_k$} \tag{ii} \\
& \framebox[23pc]{for all $k\ge0$, the vector field $k!\,\p_{k,P}$ equals
$\tpp_k$} \tag{iii}
\end{align}
\end{theorem}
\begin{proof}
We may reformulate the conditions of the theorem in terms of the
constraints $x_k$, $\alpha^u_{k,P}$ and $\alpha^v_{k,P}$ as follows:
\begin{gather}
\framebox[14pc]{for all $k\ge0$, $x_k=0$} \tag{i} \\
\framebox[14pc]{for all $k\ge0$, $\alpha_{k,P}^v=\alpha_{k,P}^u=0$} \tag{ii} \\
\framebox[14pc]{for all $k\ge0$,
$\nabla\alpha_{k,P}^v=\nabla\alpha_{k,P}^u=0$} \tag{iii}
\end{gather}

We have already shown that the Virasoro conjecture implies (i); let us
prove the converse. We argue by induction that $z_k$ vanishes, starting
with Hori's equation $z_0$. If $z_{k-1}=0$, then Theorem~\ref{master}
together with (i) implies that $\nabla z_k=0$. We conclude by \eqref{Lz}
and Proposition~\ref{Basic} that $z_k=0$.

The proof of the following lemma is analogous to that of Lemma \ref{y}.
\begin{lemma} \label{x}
We have $x_{k-1}= v\alpha_{k-1,P}^v + \DELTA\alpha_{k-1,P}^u -
k\alpha_{k,P}^v - 2\alpha_{k,Q}^v$ and
$$
\p_Q x_{k-1} = \CC_u\alpha_{k-1,P}^v + v\nabla\alpha_{k-1,P}^u -
k\nabla\alpha_{k,P}^u - 2\nabla\alpha_{k,Q}^v .
$$
\end{lemma}

By Lemma \ref{Lb} and Proposition \ref{Basic}, we see that if
$\nabla\alpha_{i,P}^*=0$ for $i\le k$, then $\alpha_{k,P}^*=0$; in
particular, (iii) implies (ii). Lemma \ref{x} shows that (ii)
implies (i) and that (i) implies (iii).
\end{proof}

\begin{corollary}
The Toda and Virasoro conjectures determine the Gromov-Witten potential $F$
of $\CP^1$ up to a constant (that is, an element of $\Q_{\eps,q}$).
\end{corollary}
\begin{proof}
The Toda and Virasoro conjectures determine the vector fields $\p_{k,Q}$
and $\p_{k,P}$ in the coordinate system $\{\p^nu,\p^nv\}_{n\ge0}$ on the
large phase space. It follows that the coordinates $\{s_k,t_k\}_{k\ge0}$
are determined up to constants of integration; but these constants are
fixed by Theorem \ref{Jet}. By inversion, we see that $u$ is determined as
a function of $\{s_k,t_k\}_{k\ge0}$. Integrating twice, using Lemma
\ref{Basic} and the string equation, we see that $F$ is determined up to a
constant.
\end{proof}

In order to determine the constant term of $F$, we may use the divisor
equation
$$
q \frac{\p F}{\p q} + \sum_{k=0}^\infty s_{k+1} \frac{\p F}{\p t_k} =
\p_QF .
$$
This fixes the constant up to an element of $\Q_\eps$; however, such
constant terms correspond to Gromov-Witten invariants of degree $0$, and no
moduli space $\Mbar_{g,n,d}$ has vanishing virtual dimension if
$d=0$. Thus, this constant vanishes.

\section{Propagating the Toda conjecture}

Consider the submanifold $\L\subset\M$ of the large phase space on
which $s_k=0$, $k>1$. Okounkov and Pandharipande \cite{OPnew} have
proved the Toda conjecture on this submanifold; that is, they prove
that $k!\,\p_{k-1,Q}$ equals $\delta_k$ along $\L$, for all $k>0$. Our
results allow us to prove that this, in conjunction with the Virasoro
conjecture, implies the full Toda conjecture.

Suppose that the constraints $\nabla\alpha_{k,Q}^v$ and
$\nabla\alpha_{k,Q}^u$ vanish to order $N$ along $\L$; the theorem of
Okounkov and Pandharipande is the case $N=1$. We now argue by
induction. The proof of Theorem \ref{Y} shows that the constraints $y_k$
vanish to order $N$ along $\L$. (Here, we use the fact that the vector
fields $\CL_{-1}$ and $\CL_0$ are tangential to $\L$, as may be seen by
inspection of the explicit formulas \eqref{L-1} and \eqref{L0}.)

Applying Theorem \ref{master}, we see that the constraints $x_k$ vanish to
order $N$ along $\L$. The proof of Theorem \ref{X} shows that the
constraints $\nabla\alpha_{k,P}^v$ and $\nabla\alpha_{k,P}^u$ vanish to
order $N$ along $\L$, in other words, that the vector fields $k!\,\p_{k,P}$
and $\tpp_k$ are equal to order $N$ along $\L$.

To prove the induction step, we must show that the vector field
$k!\,\p_{k-1,Q}-\delta_k$ vanishes to order $N+1$ along $\L$, in other
words, that $[\ell!\,\p_{\ell,P},k!\,\p_{k-1,Q}-\delta_k]$ vanishes to
order $N$ along $\L$ for all $\ell>1$. We have
\begin{multline*}
[\ell!\,\p_{\ell,P},k!\,\p_{k-1,Q}-\delta_k] =
[\ell!\p_{\ell,P},k!\,\p_{k-1,Q}] - [\tpp_\ell,\delta_k] \\ -
[\ell!\,\p_{\ell,P}-\tpp_\ell,k!\,\p_{k-1,Q}] +
[\ell!\,\p_{\ell,P}-\tpp_\ell,k!\,\p_{k-1,Q}-\delta_k] .
\end{multline*}
Obviously, $[\ell!\p_{\ell,P},k!\,\p_{k-1,Q}]$ vanishes; the commutator
$[\tpp_\ell,\delta_k]$ vanishes by Corollary \ref{magri}; the vector field
$[\ell!\,\p_{\ell,P}-\tpp_\ell,\delta_k]$ vanishes to order $N$ along $\L$,
while the vector field
$[\ell!\,\p_{\ell,P}-\tpp_\ell,k!\,\p_{k-1,Q}-\delta_k]$ vanishes to order
$2N-1\ge N$ along $\L$.

\section*{Acknowledgments}

I thank B. Dubrovin, T. Eguchi, B. Feigin, A. Orlov, R. Pandharipande,
T.~Shiota, C.-S.~Xiong, Y. Zhang and the referee for stimulating my
interest in this subject and for their helpful suggestions.

I wish to thank Kyoji Saito and Masa-Hiko Saito, and all of the other
organizers and participants in the memorable year 1999--2000 at RIMS, Kyoto
University devoted to ``Geometry of String Theory.''

The research of the author is supported in part by NSF grants DMS-9704320
and DMS-0072508.


\begin{thebibliography}{99}

\bibitem{Dorfman} I. Dorfman, ``Dirac structures and integrability of
nonlinear evolution equations.'' John Wiley, Chichester, 1993.

\bibitem{Dubrovin} B. Dubrovin, \emph{Geometry of 2D topological field
theories,} in ``Integrable systems and quantum groups, Montecalini Terme,
1993,'' eds.\ M. Francaviglia and S. Greco, Lect.\ Notes Math., vol.~1620,
Springer-Verlag, Berlin, 1996, pp.\ 120--348. \texttt{<hep-th/9407018>}

\bibitem{DZ} B. Dubrovin and Y. Zhang, in preparation.

\bibitem{EY} T. Eguchi and S.-K. Yang, \emph{The topological ${C}{\rm P}\sp
1$ model and the large-$N$ matrix integral.} Modern Phys. Lett. \textbf{A
9} (1994), 2893--2902. \texttt{<hep-th/9407134>}

\bibitem{EHY} T. Eguchi, K. Hori and S.-K. Yang, \emph{Topological $\sigma$
models and large-$N$ matrix integral.} Internat. J. Modern Phys. \textbf{A
10} (1995), 4203--4224. \texttt{<hep-th/9503017>}

\bibitem{FP} C. Faber and R. Pandharipande, \emph{Logarithmic series and
Hodge integrals in the tautological ring. With an appendix by D. Zagier.}
Michigan Math. J. \textbf{48} (2000), 215--252. \texttt{<math/0002112>}

\bibitem{FS} D.B. Fairlie and I.A.B. Strachan, \emph{The algebraic and
Hamiltonian structure of the dispersionless Benney and Toda
hierarchies}, Inverse Problems \textbf{12} (1996), 885--908.
\texttt{<math/9606022>}

\bibitem{GD} I. M. Gelfand and L. A. Dikii, \emph{Asymptotic properties of
the resolvent of Sturm-Liouville equations, and the algebra of Korteweg-de
Vries equations.} (Russian) Uspehi Mat. Nauk \textbf{30} (1975), no. 5,
67--100.  (English translation: Russian Math. Surveys \textbf{30} (1975),
no. 5, 77--113.)

\bibitem{GDorfman} I. M. Gelfand and I. Ja.\ Dorfman, \emph{Schouten bracket
and Hamiltonian operators.} (Russian) Funktsional. Anal. i
Prilozhen. \textbf{14} (1980), 71--74.

\bibitem{warsaw} E. Getzler, \emph{The Virasoro conjecture for
Gromov-Witten invariants,} ``Algebraic geometry: Hirzebruch 70 (Warsaw,
1998),'' Contemp. Math. \textbf{241}, Amer. Math. Soc., Providence, RI,
1999, pp.  147--176. \texttt{<math/9812026>}

\bibitem{darboux} E. Getzler, \emph{A Darboux theorem for Hamiltonian
operators in the formal calculus of variations,} To appear, Duke J. Math.
\texttt{<math/0002164>}

\bibitem{GP} E. Getzler and R. Pandharipande, \emph{Virasoro
constraints and the Chern classes of the Hodge bundle,} Nucl.\
Phys. \textbf{B 530} (1998), 701--714. \texttt{<math/9805114>}

\bibitem{Givental} A. Givental, \emph{Gromov - Witten invariants and
quantization of quadratic hamiltonian,} to appear, Moscow Mathematical
Journal. \texttt{<math/0108100>}

\bibitem{Hori} K. Hori, \emph{Constraints for topological strings in
$D\ge1$,} Nucl.\ Phys. \textbf{B439} (1995) 395--420.
\texttt{<hep-th/9411135>}

\bibitem{IZ} C. Itzykson and J.-B. Zuber, \emph{Combinatorics of the
modular group. II. The Kontsevich integrals,} Internat. J. Modern
Phys. \textbf{A 7} (1992), 5661--5705. \texttt{<hept-th/9201001>}

\bibitem{Kontsevich} M. Kontsevich, \emph{Intersection theory on moduli
spaces of curves and the matrix Airy function},
Commun. Math. Phys. \textbf{147} (1992), 1--23.

\bibitem{Kon} M. Kontsevich, \emph{Enumeration of rational curves via torus
actions.} In ``The moduli space of curves (Texel Island, 1994),'' 335--368,
Progr. Math. \textbf{129}, Birkh\"auser Boston, Boston, MA, 1995.

\bibitem{K} B. A. Kupershmidt, \emph{Discrete Lax equations and
differential-difference calculus,} Ast\'erisque \textbf{123} (1985).

\bibitem{L} E. Looijenga, \emph{Intersection theory on Deligne-Mumford
compactifications (after Witten and Kontsevich).} S\'eminaire Bourbaki,
Vol.\ 1992/93. Ast\'erisque No.\ 216 (1993), Exp.\ No.\ 768, pp.\ 187--212.

\bibitem{Manin} Y. I. Manin, \emph{Algebraic aspects of nonlinear
differential equations.} (Russian) In ``Current problems in mathematics,''
Vol. 11, pp. 5--152. Akad. Nauk SSSR Vsesojuz. Inst. Nau\v cn. i
Tehn. Informacii, Moscow, 1978.

\bibitem{frobenius} Y. I. Manin, ``Frobenius manifolds, quantum cohomology,
and moduli spaces.'' American Mathematical Society Colloquium Publications
\textbf{47}. American Mathematical Society, Providence, RI, 1999.

\bibitem{OP} A. Okounkov and R. Pandharipande, \emph{Gromov-Witten theory,
Hurwitz numbers, and matrix models, I.} \texttt{<math/0101147>}

\bibitem{OPnew} A. Okounkov and R. Pandharipande, private communication.

\bibitem{Olver} P. Olver, ``Applications of Lie groups to differential
equations.'' Graduate Texts in Mathematics, v.\ 107. Springer-Verlag, New
York, 1993.

\bibitem{P} R. Pandharipande, \emph{The Toda equations and the
Gromov-Witten theory of the Riemann sphere.} Lett. Math. Phys. \textbf{53}
(2000), 59--74. \texttt{<math/9912166>}

\bibitem{TT} K. Takasaki and T. Takebe, \emph{Quasi-classical limit of Toda
hierarchy and $W$-infinity symmetries.}  Lett. Math. Phys. \textbf{28}
(1993), 165--176.

\bibitem{Witten} E. Witten, \emph{Two dimensional gravity and
intersection theory on moduli space}, Surveys in Differential
Geom. \textbf{1} (1991), 243--310.

\end{thebibliography}
\end{document}